\documentclass[11pt]{amsart}
\usepackage{url}
\usepackage{hyperref,amssymb}
\hypersetup{colorlinks=true,linkcolor=blue,
urlcolor=cyan,citecolor=magenta} 
\usepackage[hyperpageref]{backref}
\usepackage{frcursive,calrsfs,mathrsfs}
\usepackage{xcolor}
\newtheorem{theorem}{Theorem}[section]

\newtheorem{conjecture}[theorem]{Conjecture}

\numberwithin{equation}{section}
\newcommand{\Q}{{\mathbb{Q}}}
\newcommand{\Z}{{\mathbb{Z}}}

\newcommand{\ds}{\displaystyle}

\newcommand{\wh}{\widehat}
\newcommand{\ft}{\footnotesize}

\newcommand{\CH}{{\mathcal H}}
\newcommand{\CE}{{\mathcal E}}
\newcommand{\CF}{{\mathcal F}}
\newcommand{\CA}{{\mathcal A}}

\newcommand{\CX}{{\mathcal X}}
\newcommand{\CC}{{\mathcal C}}
\newcommand{\CG}{{\mathcal G}}
\newcommand{\order}{\raise0.8pt \hbox{${\scriptstyle \#}$}}

\newcommand{\lien}{\mathrel{\mkern-4mu}}
\newcommand{\too}{\relbar\lien\rightarrow}
\newcommand{\tooo}{\relbar\lien\relbar\lien\too}

\newcommand{\plus}{\ds\mathop{\raise 0.5pt \hbox{$\bigoplus$}}\limits}
\newcommand{\prd}{\ds\mathop{\raise 1.0pt \hbox{$\prod$}}\limits}
\newcommand{\sm}{\ds\mathop{\raise 1.0pt \hbox{$\sum$}}\limits}
\newcommand{\ffrac}[2]{\hbox{\ft $\displaystyle\frac{#1}{#2}$}}
\newcommand{\Gal}{{\rm Gal}}
\newcommand{\Norm}{\hbox{\bf$\textbf{N}$}}
\newcommand{\BJ}{\hbox{\bf$\textbf{J}$}}

\newcommand{\Ker}{{\rm Ker}}
\newcommand{\ram}{{\rm ram}}
\newcommand{\ab}{{\rm ab}}
\newcommand{\ar}{{\rm ar}}
\newcommand{\Kappa}{{\hbox{\large $\kappa$}}}
\newcommand{\BKappa}{\hbox{\Large $\kappa$}}
\newcommand{\Nu}{\hbox{\Large $\nu$}}

\author[Georges Gras]{Georges Gras}

\address{\hspace{-0.5cm}Villa la Gardette, 4 chemin Ch\^ateau Gagni\`ere, 
38520 Le Bourg d'Oisans \rm {\url{http://orcid.org/0000-0002-1318-4414} } }
\email{g.mn.gras@wanadoo.fr}



\begin{document}

\title[Exceptional case of the Main Conjecture]
{Exceptional case of the non semi-simple \\ 
 Real Abelian Main Conjecture}

\date{August 8, 2023}

\begin{abstract} In the papers: ``The Chevalley--Herbrand formula and the real 
abelian Main Conjecture (New criterion using capitulation of the class group),
\textit{J. Number Theory} \textbf{248} (2023)'' and ``On the real abelian main 
conjecture in the non semi-simple case, arXiv (2023)'', we consider real cyclic 
fields $K$ and primes $\ell \equiv 1 \pmod {2p^N}$ totally inert in $K$, implying
implicitly, $K \cap \Q(\mu_{p^\infty}^{})^+ = \Q$. In this complementary work, we 
examine the non linearly disjoint case.
\end{abstract}

\maketitle

\vspace{-0.6cm}

\tableofcontents

\vspace{-0.6cm}
\section{Introduction}

Let $K/\Q$ be a real cyclic extension of Galois group $\CG_K \simeq \gamma 
\oplus \Gamma$, with $\gamma$ of prime-to-$p$ order $d$ and $\Gamma 
\simeq \Z/p^c \Z$, $c \geq 1$. Let $\CH_K$ be the $p$-class group of $K$ as 
$\Z_p[\CG_K]$-module, and let $E_K$ (resp. $F_K$) be the group of units
(resp. of Leopoldt's cyclotomic units) of $K$; set $\CE_K = E_K \otimes \Z_p$,
$\CF_{\!K} = F_K \otimes \Z_p$. We consider $\wh \CE_K = \wh E_K \otimes \Z_p$,
where $\wh E_K$ is the subgroup of $E_K$ generated by the units of the
strict subfields of $K$.

\smallskip
The ``Real Abelian Main Conjecture'' (RAMC for short), stated in \cite{Gra0}, 
writes as follows in terms of arithmetic $\varphi$-objects 
(for irreducible $p$-adic character $\varphi  = \varphi_0^{} \varphi_p$ of $K$,
$\varphi_0^{}$ of prime-to-$p$ order, $\varphi_p$ of $p$-power order):
\begin{equation}\label{RAMC}
\order \CH^\ar_{K,\varphi} = \order (\CE_K / \wh \CE_K \, \CF_{\!K})_{\varphi_0^{}} ,
\end{equation}

\noindent
where $\CH^\ar_{K,\varphi} :=  \{x \in \CH_K,\ \,\Norm_{K/\Kappa}(x) = 1,\, 
\forall \, \BKappa \varsubsetneqq K \}_{\varphi_0^{}}$ and
where the semi-simple ${\varphi_0^{}}$-components are given by the idempotents
$e_{\varphi_0^{}} := \ffrac{1}{\order \gamma} \, \sm_{\tau \in \gamma} 
\varphi_0^{}(\tau^{-1}) \tau$.

\smallskip
We have proved in \cite{Gra1,Gra2} that the RAMC \eqref{RAMC} holds for $K$ 
as soon as there exists a prime number $\ell \equiv 1\! \pmod {2p^N}$, totally inert 
in $K$, such that $\CH_K$ capitulates in the auxiliary $p$-sub-extension $L/K$ of 
$K(\mu_\ell^{})/K$. 
Obviously, when $\CC := K \cap \Q(\mu_{p^\infty}^{})^+$
is distinct from $\Q$, any $\ell \equiv 1\! \pmod {2p^N}$ totally splits in $\CC$ 
and can not be totally inert in $K$ (this was not made explicit in our papers);
so we shall look at this particular case using a similar process of capitulation
of the $p$-class group $\CH_K$ in suitable absolute abelian extensions.

\smallskip
One proves the RAMC, except in the case $K \subset \Q(\mu_{p^\infty}^{})^+$ 
using Greenberg's conjecture, and in the ``special case'' $K \cap
\Q(\mu_{p^\infty}^{})^+ \varsubsetneqq K$ with $\varphi_0^{}$ 
character of $\Q(\mu_p)^+$, which remains open.

\smallskip
Some proofs, identical to that given in the two previous papers, will not be
given again, as well as notations, definitions and references.

\section{Recalls about real abelian fields}
Let $\Q^\ab$ be the maximal abelian extension of $\Q$. We consider the set of 
even irreducible $p$-adic characters $\Phi$, for a prime $p \geq 2$, and the set 
$\CX$ of irreducible rational characters of $\Q^\ab$; then the subsets 
of irreducible characters $\Phi_K$, $\CX_K$, of $K \subset \Q^\ab$. Thus, $K$ is 
the fixed field of a unique $\chi \in \CX$ and we put $K = K_\chi$.

\subsection{Arithmetic \texorpdfstring{$\chi$ and $\varphi$-objects}{Lg}}\label{characters}
Let $k$ be the subfield of $K = K_\chi$ fixed by $\Gamma \simeq \Z/p^c\Z$.
The field of values of $\psi \mid \chi$, of degree~$1$, is $\Q(\mu_{dp^c})$, 
direct compositum $\Q(\mu_d) \Q(\mu_{p^c})$; thus $\psi = \psi_0 \psi_p$, 
$\psi_0$ of order $d$,  $\psi_p$ of order $p^c$ and $\chi = \chi_0^{} \chi_p$, 
$\chi_0^{} \in \CX_k$ above $\psi_0$ and $\chi_p \in \CX_{K^\gamma}$ 
above $\psi_p$. Similarly, in the direct compositum $\Q_p(\mu_d) \Q_p(\mu_{p^c})$,
the irreducible $p$-adic characters $\varphi \mid \chi$ are of the form
$\varphi = \varphi_0^{} \varphi_p$, $\varphi_0^{} \mid \chi_0^{}$, and
where $\varphi_p = \chi_p$ since $\Gal(\Q_p(\mu_{p^c})/\Q_p)\simeq 
\Gal(\Q(\mu_{p^c})/\Q)$. 

\smallskip
If $\CA$ is a $\Z_p[\CG_K]$-module, $\CA_{\varphi_0^{}} := \CA^{e_{\varphi_0^{}}}$ 
denotes its semi-simple ${\varphi_0^{}}$-component given by the idempotent
$e_{\varphi_0^{}} := \ffrac{1}{\order \gamma} \, \sm_{\tau \in \gamma} 
\varphi_0^{}(\tau^{-1}) \tau \in \Z_p[\gamma]$.

Recall that for $K = K_\chi$ and $\varphi \mid \chi$:
\begin{equation}\label{maindef}
\left \{\begin{aligned}
\CH^\ar_{K,\chi} & :=  \{x \in \CH_K,\ \,\Norm_{K/\Kappa}(x) = 1,\, 
\forall \, \BKappa \varsubsetneqq K \}, \\
\CH^\ar_{K,\varphi} & :=  (\CH^\ar_{K,\chi})_{\varphi_0^{}}.
\end{aligned}\right.
\end{equation}

The above definitions of $\chi$ and $\varphi$-objects lead to the fundamental 
formula \cite[Corollary 2.3]{Gra2}:

\begin{theorem}\label{isotopicphi}
For all $\rho \in \CX$, we have the semi-simple decomposition
$\CH^\ar_{K_\rho,\rho} = \plus_{\varphi \mid \rho} \CH^\ar_{K_\rho,\varphi}$, 
and for all real cyclic extension $K/\Q$, we have the global formula
$\order \CH_K = \prd_{\rho \in \CX_K} \order \CH^\ar_{K_\rho,\rho} =
\prd_{\rho \in \CX_K}\prd_{\varphi \mid \rho}\order \CH^\ar_{K_\rho,\varphi}$.
\end{theorem}

\subsection{Analytic expression of 
\texorpdfstring{$\order \CH^\ar_{K,\chi}$}{Lg} and the RAMC}\label{subIII3}

Since the RAMC is trivial for $\chi$ of $p$-power order, 
we assume $d > 1$, $\chi_0^{} \ne 1$. 
The interpretation of Leopoldt's analytic formulas yields the following result
for $K=K_\chi$ real  (see the survey \cite[Theorem 3.2]{Gra2}):
\begin{theorem}\label{chiformulaH}
Let $\CH^\ar_{K,\chi} := \{x \in \CH_K,\,  
\Norm_{K/\Kappa}(x) = 1,\, \hbox{$\forall \ \BKappa \varsubsetneqq K$} \}$. 
Then $\order \CH^\ar_{K,\chi} =  \big (\CE_K\! : \wh \CE_K \, \CF_{\!K} \big)$ where
$\wh \CE_K$ is the subgroup of $\CE_K$ generated by the $\CE_{\Kappa}$'s for 
all $\BKappa \varsubsetneqq K$ and $\CF_{\!K}$ is the group of Leopoldt's
cyclotomic units.
\end{theorem}

We have formulate in the 1970's \cite{Gra0} the following RAMC:
\begin{conjecture}\label{mainconj}
Let $\chi = \chi_0^{} \chi_p$, $\chi_0^{} \ne 1$, be an even irreducible rational character 
and let $K=K_\chi$. For all irreducible $p$-adic character $\varphi \mid \chi$
(whence $\varphi_0^{} \mid \chi_0^{}$, $\varphi_p = \chi_p$), we have
$\order \CH^\ar_{K,\varphi} = \order (\CE_K / \wh \CE_K \, \CF_{\!K})_{\varphi_0^{}}$
(${\varphi_0^{}}$-component of the $\chi$-object $\CE_K / \wh \CE_K \, \CF_{\!K}$). 
\end{conjecture}

As we have explained, the particular case $K \cap \Q(\mu_{p^\infty}^{})^+ \ne \Q$
must be examined separately; this is coherent with the specific properties of
$p$-adic ${\rm L}$-functions regarding characters conductors powers of $p$.

\smallskip
Let $\CC := K \cap \Q(\mu_{p^\infty}^{})^+$ and $\CC_0 := k \cap \Q(\mu_p^{})^+$
with $d = [k : \Q] > 1$, ${\varphi_0^{}} \mid \chi_0^{} \ne 1$ and $[K : k] = p^c > 1$. 
Note that for $p \leq 3$, $\CC_0 = \Q$ and $\CC$ is a layer of the cyclotomic 
$\Z_2$-extension $\Q_\infty$; since $k \ne \Q$ and ${\varphi_0^{}} \ne 1$, this 
case will be solved in \S\,\ref{known}. So we can focus on the case $p > 3$ 
in the sequel.

\smallskip
We will analyze the RAMC, under capitulation properties, in two steps
corresponding to the following cases:

\smallskip
(i) $\CC = K \subset \Q(\mu_{p^\infty}^{})^+$ (cf. Section \ref{green}),

\smallskip
(ii) $\Q \varsubsetneqq \CC \varsubsetneqq K$, with
the two sub-cases $\CC_0 \varsubsetneqq k$ and the ``special case'' $\CC_0 = k$
(cf. Section \ref{special}).

\section{Case \texorpdfstring{$K \subset \Q(\mu_{p^\infty}^{})^+$}{Lg} 
and Greenberg's conjecture}\label{green}
Note that this framework makes sense only if Vandiver's conjecture is false,
otherwise all the subfields of $\Q(\mu_{p^\infty}^{})^+$ are $p$-principal.

\smallskip
Thus, $k$ is a subfield of  $\Q(\mu_p^{})^+$ and $K$ is the layer $k_c$ of the 
cyclotomic $\Z_p$-extension $k_\infty$ of $k$.

\unitlength=0.52cm
$$\vbox{\hbox{\hspace{0.6cm} 
\begin{picture}(11.5,2.6)
\put(7.5,2.4){\line(1,0){3.6}}
\put(0.7,2.4){\line(1,0){3.6}}
\bezier{150}(-1.6,1.0)(-1.0,1.55)(-0.4,2.1)
\put(-1.6,0.5){\line(1,0){4.15}}
\put(-3.6,-1.5){\line(1,0){3.4}}
\bezier{150}(1.3,-1.0)(1.9,-0.45)(2.5,0.1)
\bezier{150}(-3.6,-1.0)(-3.0,-0.45)(-2.4,0.1)
\bezier{150}(8.3,-1.0)(8.9,-0.45)(9.5,0.1)
\put(4.6,0.50){\line(1,0){4.8}}
\put(2.45,-1.50){\line(1,0){4.8}}
\bezier{150}(3.3,1.)(3.9,1.55)(4.5,2.1)
\bezier{150}(10.3,1.)(10.9,1.55)(11.5,2.1)
\put(2.7,0.35){\ft$K\!=\!k_c$}
\put(-4.3,-1.6){\ft$\Q$}
\put(0.1,-1.6){\ft$K^\gamma\!=\! \Q_c$}
\put(7.6,-1.6){\ft$\Q_n$}
\put(8.4,-1.50){\line(1,0){3.3}}
\put(11.9,-1.6){\ft$\Q_\infty$}
\put(10.4,0.45){\line(1,0){2.6}}
\put(13.1,0.35){\ft$k_{\infty}$}
\put(-2.2,0.35){\ft$k$}
\put(0.2,0.65){\ft$\Gamma$}
\put(2.15,-0.55){\ft$\gamma$}
\bezier{250}(-3.6,-1.3)(0.0,-0.75)(2.5,0.35)
\put(-0.8,-0.5){\ft$\CG_K$}%
\put(4.6,2.3){\ft$\Q(\mu_{p^{c+1}})^+$}
\put(9.6,0.35){\ft$k_n$}
\put(11.5,2.3){\ft$\Q(\mu_{p^{n+1}})^+$}
\put(-1.1,2.3){\ft$\Q(\mu_p)^+$}
\end{picture}   }} $$
\vspace{0.6cm}
\unitlength=1.0cm

\noindent
This case can not use an auxiliary totally inert prime $\ell \equiv 1 \pmod {2p^N}$ 
(since it is totally split in a non-trivial subfield of $K$), but Greenberg's conjecture 
implying a capitulation principle in $k_\infty$ (in some sense, $p$ replaces $\ell$ 
and $k_\infty$ replaces $L \subset K(\mu_\ell^{})$). 

\smallskip
Assuming Greenberg's conjecture \cite[Theorem 2]{Gre}, we know that, for all 
$m \geq n \geq n_0$, the arithmetic norms $\Norm_{k_m/k_n} : \CH_{k_m} \to \CH_{k_n}$ 
are isomorphisms and that $\CH_{k_m} = \CH_{k_m}^{\Gal(k_m/k_n)}$;
so, $\Nu$ denoting the algebraic norm,
$\Nu_{k_m/k_n}(\CH_{k_m}) = \BJ_{k_m/k_n} \circ \Norm_{k_m/k_n}(\CH_{k_m})
= \CH_{k_m}^{p^{m-n}} = \BJ_{k_m/k_n}(\CH_{k_n}) = 1$ for $m \gg n \geq n_0$.
It follows that the $p$-class group of any layer $k_n$ 
($n \geq 0$) capitulates in $k_\infty$; so, $\CH_K$ capitulates
in $k_\infty$. 

\smallskip
Since $p$ is totally ramified in $k_\infty/\Q$, the Chevalley--Herbrand 
formula \cite[pp. 402-406]{Che} yields $\order \CH_{k_n}^{\Gal(k_n/k)} 
= \order \CH_{k}$, and in the same way as for \cite[Corollary 4.5]{Gra2},
we have: 
$$\CH_{k_n,\varphi_0^{}}^{\Gal(k_n/k)} \simeq \big( \CE_{k}/
\Norm_{k_n/k}(\CE_{k_n}) \big)_{\varphi_0^{}}, $$
of order $\order \CH_{k,\varphi_0^{}}$, as soon as $\CH_{k,\varphi_0^{}}$
capitulates in $k_n$.

\smallskip
In terms of cyclotomic units $\CF_{k_n}$, we have the fundamental
norm relation $\Norm_{k_m/k_n}(\CF_{k_m}) = \CF_{k_n}$, for all
$m \geq n \geq 0$; the principle of computation giving \cite[Theorem 4.4]{Gra2} 
holds, which leads to the proof of the non semi-simple RAMC \ref{mainconj}, 
for $K \subset \Q(\mu_{p^\infty}^{})^+$.

\smallskip
However, Greenberg's conjecture is far to be proved and is probably more difficult 
than the capitulation conjecture using primes $\ell$ totally inert in $K/\Q$, since 
these are infinite in number and many numerical computations go in this direction.

\section{Analysis of the RAMC when \texorpdfstring{$K \cap \Q(\mu_{p^\infty}^{})^+ 
\varsubsetneqq K$}{Lg} }

\subsection{Expression of \texorpdfstring{$\order \CH^\ar_{K,\varphi}$}{Lg}} \label{special}

We will use the following definitions when $K = K_\chi$ (of degree 
$d p^c$) contains $\CC := K \cap \Q(\mu_{p^\infty}^{})^+$, such that 
$\Q \varsubsetneqq \CC \varsubsetneqq K$. Let $\varphi = \varphi_0 \varphi_p \mid \chi$.

\smallskip
Let $M_0 := M(\ell) \CC$, where $M(\ell)$ is the subfield of degree $p^N$ of 
$\Q(\mu_\ell^{})$ for $\ell \equiv 1 \pmod {2p^N}$, $N \geq 1$, with $\ell$ prime 
totally inert in $K/\CC$ (possible from the density theorem) and let $L= M_0 K$.
So any prime ideal ${\mathfrak l} \mid \ell$ of $K$ is totally ramified in $L/K$, 
then (totally if $N$ is large enough) split in $\CC$ and totally inert in $L/M_0$:
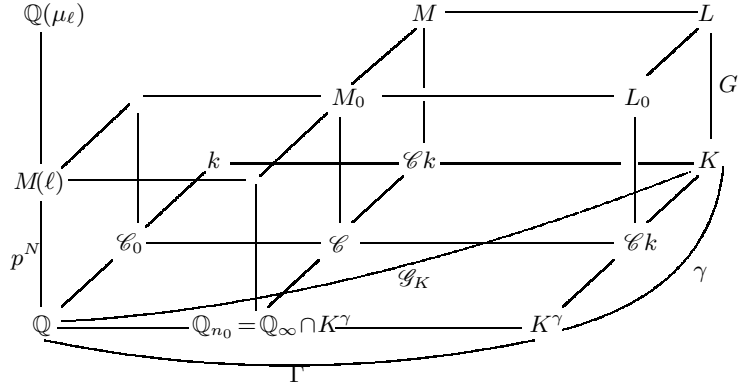
\unitlength=0.56cm
$$\vbox{\hbox{\hspace{2.4cm} 
\begin{picture}(11.5,6.2)
\put(5.5,6){\line(1,0){6}}
\put(4.0,4){\line(1,0){5.5}}
\put(5.4,2.4){\line(1,0){4.25}}
\put(10.1,2.4){\line(1,0){1.3}}
\put(0.4,2.4){\line(1,0){4.0}}
\bezier{150}(-1.6,0.9)(-0.9,1.55)(-0.2,2.2)
\put(-1.6,0.5){\line(1,0){4.15}}
\put(-3.7,-1.5){\line(1,0){3.1}}
\bezier{150}(1.3,-1.0)(1.9,-0.45)(2.5,0.1)
\bezier{150}(-3.6,-1.0)(-3.0,-0.45)(-2.4,0.1)
\put(-1.65,4.0){\line(1,0){4.2}}
\put(-3.45,2.0){\line(1,0){4.2}}
\bezier{150}(1.1,2.1)(1.9,2.85)(2.7,3.6)
\bezier{150}(-3.6,2.3)(-2.8,3.05)(-2.0,3.8)
\bezier{150}(8.3,-1.0)(8.9,-0.45)(9.5,0.1)
\put(3.5,0.50){\line(1,0){6}}
\put(2.9,-1.50){\line(1,0){4.4}}
\bezier{150}(3.3,1.)(3.9,1.55)(4.5,2.1)
\bezier{150}(3.3,4.3)(4.05,4.95)(4.8,5.6)
\bezier{150}(10.3,1.)(10.9,1.55)(11.5,2.1)
\bezier{150}(10.3,4.5)(10.9,5.05)(11.5,5.6)
\put(3.0,1.0){\line(0,1){2.5}}
\put(1.0,-1.1){\line(0,1){3.0}}
\put(-4.1,-1.0){\line(0,1){2.5}}
\put(-1.8,1.0){\line(0,1){2.6}}
\put(-4.75,1.8){\ft$M\!(\ell)$}
\put(-4.1,2.4){\line(0,1){3.1}}
\put(-4.5,5.8){\ft$\Q(\mu_\ell)$}
\put(5.0,2.9){\line(0,1){0.98}}
\put(5.0,4.2){\line(0,1){1.4}}
\put(10.0,1.0){\line(0,1){2.5}}
\put(11.8,3.0){\line(0,1){2.5}}
\put(2.75,0.3){\ft$\CC$}
\put(-4.3,-1.6){\ft$\Q$}
\put(-0.5,-1.6){\ft$\Q_{n_0}\!=\!\Q_\infty\! \cap\! K^\gamma$}
\put(7.5,-1.6){\ft$K^\gamma$}
\put(-0.15,2.3){\ft$k$}
\put(4.5,2.3){\ft$\CC k$}
\put(-2.3,0.35){\ft$\CC_0$}
\put(9.75,0.35){\ft$\CC k$}
\put(11.5,2.3){\ft$K$}
\put(-4.8,0.0){\ft$p^N$}
\put(2.8,3.8){\ft$M_0$}
\put(4.7,5.8){\ft$M$}
\put(9.75,3.8){\ft$L_0$}
\put(11.5,5.8){\ft$L$}
\put(12.0,4.1){\ft$G$}
\put(4.4,-0.5){\ft$\CG_K$}
\bezier{300}(-3.6,-1.35)(3.1,-1.0)(11.3,2.2)
\bezier{300}(-4.0,-1.8)(1.8,-3.0)(7.6,-1.8)
\put(1.8,-2.8){\ft$\Gamma$}
\bezier{200}(8.3,-1.55)(11.8,-0.6)(12.1,2.3)
\put(11.4,-0.3){\ft$\gamma$}
\end{picture}   }} $$
\vspace{1.0cm}
\unitlength=1.0cm

We use some notations of \cite{Gra2}; in particular,
let $\wh K$ and $\wh L$ be the subfields 
of $K$ and $L$ such that $[K : \wh K] = [L : \wh L] = p$.
Put $\wh \Norm := \Norm_{L/\wh L}$, $\wh \BJ := \BJ_{L/\wh L}$.
From \cite[Lemma 4.3]{Gra2}, if $\CH_K$ capitulates in $L$, then 
$\CH_{\wh K}$ capitulates in $\wh L$. By abuse, put $(A : B)_{\varphi_0^{}}
= \order (A/B)_{\varphi_0^{}}$.

\smallskip
Under the above assumptions, the $p$-localized Chevalley--Herbrand formulas 
for $p$-class groups are the following ones (to get the orders, for $\varphi_0^{} \ne 1$, 
of the $\varphi_0^{}$-components $\CH_{L,\varphi_0^{}}^G$, $G := \Gal(L/K)$, we have 
followed the process given in Jaulent's Thesis \cite[Chapitre III, p. 167]{Jau} and 
described again in \cite[Section 4.1]{Gra2}):
\begin{equation}\label{chevalley}
\left \{\begin{aligned}
\order \CH_{L,\varphi_0^{}}^G &= \ds  \order \CH_{K,\varphi_0^{}} 
\times \frac{p^{N \, r_{\varphi_0^{}}}}
{(\CE_K : \CE_K \cap \Norm_{L/K}(L^\times))_{\varphi_0^{}}} \\
\order \CH_{\wh L,\varphi_0^{}}^G &= \ds  \order \CH_{\wh K,\varphi_0^{}} 
\times \frac{p^{N \, r_{\varphi_0^{}}}}
{(\CE_{\wh K} : \CE_{\wh K} \cap \Norm_{\wh L/\wh K}(\wh L^\times))_{\varphi_0^{}}},
\end{aligned}\right.
\end{equation}

\noindent
$r_{\varphi_0^{}} = [\CC : \CC_0] \, d_{\varphi_0^{}}$, where 
$d_{\varphi_0^{}} = [\Q_p(g_{\varphi_0^{}}): \Q_p]$ as usual. 

\smallskip
Under capitulations in $L/K$ and $\wh L/\wh K$, we have the exact sequences:
\begin{equation}\label{ram}
\left \{\begin{aligned}
1 \to  \CH_{L,\varphi_0^{}}^\ram  \too &\ \CH_{L,\varphi_0^{}}^{G} \\
& \too (\CE_K\cap \Norm_{L/K}(L^\times)/\Norm_{L/K}(\CE_L))_{\varphi_0^{}} \to 1, \\
1 \to \CH_{\wh L,\varphi_0^{}}^\ram \too &\ \CH_{\wh L,\varphi_0^{}}^{G} \\
& \too (\CE_{\wh K}\cap \Norm_{\wh L/\wh K}({\wh L}^\times)/
\Norm_{\wh L/\wh K}(\CE_{\wh L}))_{\varphi_0^{}} \to 1,
\end{aligned}\right.
\end{equation}

\noindent
where $\CH^\ram$ is the subgroup of $\CH$ generated by the $p$-classes 
of the ramified prime ideals ${\mathfrak L} \mid{\mathfrak l} \mid \ell$ (for $L$ and $\wh L$).

\smallskip
The $\Z_p[\gamma]$-modules $\Omega_K = \{(s_i)_{i=1}^r \in G^r, \ \,
\prd_{i=1}^r s_i = 1\}$ and $\Omega_{\wh K} = \{(\wh s_i)_{i=1}^r \in \wh G^r,\ \,
\prd_{i=1}^r \wh s_i = 1\}$, where $r = [\CC : \Q]$, are isomorphic. Then, from 
\eqref{chevalley}, \eqref{ram} giving two expressions of $\CH_{L,\varphi_0^{}}^{G}$
and $\CH_{\wh L,\varphi_0^{}}^{G}$:
\begin{equation}\label{relations}
\left \{\begin{aligned}
\order  \CH_{K,{\varphi_0^{}}} & = \frac{\order \CH_{L,{\varphi_0^{}}}^\ram \times 
(\CE_K : \Norm (\CE_L))_{\varphi_0^{}}}{\order \Omega_{K,\varphi_0^{}}}, \\
\order  \CH_{\wh K,{\varphi_0^{}}} &= \frac{\order \CH_{\wh L,{\varphi_0^{}}}^\ram \times
(\CE_{\wh K} : \Norm (\CE_{\wh L}))_{\varphi_0^{}}}
{\order \Omega_{\wh K, \varphi_0^{}}}, 
\end{aligned}\right.
\end{equation}

One sees easily that $\CH_{L,{\varphi_0^{}}}^\ram = \wh \BJ \big(\CH_{\wh L,{\varphi_0^{}}}^\ram \big)$
since prime ideals over $\ell$ are inert in $L/\wh L$.
We have the exact sequence defining $\CH^\ar_{K,\varphi}$:
$$1 \to (\CH^\ar_{K,\chi})_{\varphi_0^{}} =  \CH^\ar_{K,\varphi} \tooo \CH_{K,{\varphi_0^{}}} 
\mathop{\tooo}^{\wh {\bf N}} \CH_{\wh K,{\varphi_0^{}}} \to 1 ; $$
then, from \eqref{relations}, and the fact that $\order \Omega_{K,{\varphi_0^{}}} 
= \order \Omega_{\wh K,{\varphi_0^{}}}$, we get, under capitulations, the fundamental relation:
\begin{equation}\label{fondamental}
\order \CH^\ar_{K,\varphi} = \frac{\order \wh \BJ\big (\CH_{\wh L,{\varphi_0^{}}}^\ram \big)}
{\order \CH_{\wh L,{\varphi_0^{}}}^\ram}
\times 
\frac{( \CE_K : \Norm (\CE_L))_{\varphi_0^{}}}{(\CE_{\wh K} : \Norm (\CE_{\wh L}))_{\varphi_0^{}}},
\end{equation}

\noindent
to be compared with $X_{{\varphi_0^{}}} := \big(\CE_K : \CE_{\wh K} \CF_{\!K}\big)_{{\varphi_0^{}}}$.
So, the proof of the RAMC is essentially based on $\Norm(\CF_L)$
regarding $\Norm (\CE_L)$.

\smallskip
We have $\Norm(\CF_L) = \CF_K^{1-\tau_\ell^{-1}}$, where $\tau_\ell^{}$ is the
Artin automorphism of $\ell$ in $K/\Q$; it generates $\Gal(K/\CC)$ and is 
of $p$-power order if and only if $\CC_0 = k$, which defines the ``special case''.

\subsection{Proof of the RAMC when \texorpdfstring{$\CC_0
\varsubsetneqq k $}{Lg}}\label{known}

In this case, $(1-\tau_\ell^{-1}) e_{\varphi_0^{}}$ is invertible in
$\Z_p[\CG_K] e_{\varphi_0^{}}$ \cite[Proposition 3.4]{Gra2}.
Recall that this applies for $p \leq 3$ since $\CC_0 = \Q$ and $k \ne \Q$.

\smallskip
The computations \cite[Proof of Theorem 4.4]{Gra2} leads in this case to the relation:
\begin{equation}\label{X}
X_{{\varphi_0^{}}} = \frac{(\CE_K : \Norm (\CE_L))_{{\varphi_0^{}}}}
{(\CE_{\wh K} : \Norm (\CE_{\wh L}))_{{\varphi_0^{}}}} \times Z_{{\varphi_0^{}}},
\end{equation}
where:
\begin{equation*}
Z_{{\varphi_0^{}}} :=  \frac{(\CE_{\wh K} \Norm (\CE_L) : \Norm (\CE_{\wh L}) \CF_{\!K})_{{\varphi_0^{}}}}
{(\CE_{\wh K} \Norm (\CE_L) : \Norm (\CE_L))_{{\varphi_0^{}}}}
\times (\Norm (\CE_{\wh L}) \wh \CF_{\!K} : \Norm (\CE_{\wh L}))_{{\varphi_0^{}}},
\end{equation*}

\noindent
with $\wh \CF_{\!K} := \CF_{\!K} \cap \CE_{\wh K}$.
Since $\Norm (\CE_{\wh L}) \CF_{\!K} \subseteq \Norm (\CE_L)$,  
$Z_{{\varphi_0^{}}}$ is an integer and we obtain 
from the relations \eqref{fondamental} and \eqref{X} under capitulations:
$$X_{{\varphi_0^{}}} = \order \CH^\ar_{K,\varphi} \times 
\frac{\order \CH_{\wh L,{\varphi_0^{}}}^\ram}
{\order \wh\BJ \big(\CH_{\wh L,{\varphi_0^{}}}^\ram\big)} \times Z_{\varphi_0^{}} 
\geq \order \CH^\ar_{K,\varphi}; $$
in other words, one gets the supplementary factor $\order \Ker (\wh\BJ)_{\varphi_0^{}}$
which was trivial in the case $\ell$ totally inert in $K$ and $\varphi_0^{} \ne 1$. Whence 
the equality as usual implying the supplementary interesting relations:

\smallskip
$\ \bullet$ $\CF_{\!K,\varphi_0^{}} \cap \CE_{\wh K,\varphi_0^{}} \subseteq \Norm (\CE_{\wh L,\varphi_0^{}})$,

\smallskip
$\ \bullet$ $\Norm (\CE_{\wh L,\varphi_0^{}}) \CF_{\!K,\varphi_0^{}} = \Norm (\CE_{L,\varphi_0^{}})$,

\smallskip
$\ \bullet$ $\wh \BJ \big(\CH_{\wh L,\varphi_0^{}}^\ram \big) \simeq  \CH_{\wh L,\varphi_0^{}}^\ram$.

\subsection{The special case \texorpdfstring{$\CC_0 = k$}{Lg}} \label{known2}
So, $\tau_\ell^{}$ generates $\Gal(K/\CC)$, of $p$-power order, and 
$(1-\tau_\ell^{-1}) e_{\varphi_0^{}}$ is not invertible in $\Z_p[\CG_K] e_{\varphi_0^{}}$.
\unitlength=0.62cm
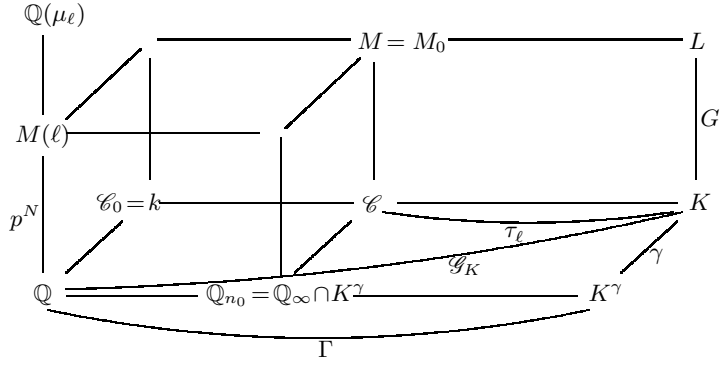
$$\vbox{\hbox{\hspace{3.4cm} 
\begin{picture}(11.5,5.0)
\put(4.6,4){\line(1,0){5.0}}
\put(-1.6,0.5){\line(1,0){4.15}}
\put(-3.6,-1.5){\line(1,0){2.8}}
\bezier{150}(1.3,-1.0)(1.9,-0.45)(2.5,0.1)
\bezier{150}(-3.6,-1.0)(-3.0,-0.45)(-2.4,0.1)
\put(-1.65,4.0){\line(1,0){4.15}}
\put(-3.6,2.0){\line(1,0){4.15}}
\bezier{150}(1.1,2.1)(1.9,2.85)(2.7,3.6)
\bezier{150}(-3.6,2.3)(-2.8,3.05)(-2.0,3.8)
\bezier{150}(8.3,-1.0)(8.9,-0.45)(9.5,0.1)
\put(3.5,0.50){\line(1,0){6}}
\put(2.55,-1.50){\line(1,0){4.8}}
\put(3.0,1.0){\line(0,1){2.5}}
\put(1.0,-1.1){\line(0,1){3.0}}
\put(9.9,1.0){\line(0,1){2.6}}
\put(-1.8,1.0){\line(0,1){2.6}}
\put(-4.1,-1.0){\line(0,1){2.5}}
\put(-4.7,1.8){\ft$M(\ell)$}
\put(-4.1,2.4){\line(0,1){1.7}}
\put(-4.5,4.4){\ft$\Q(\mu_\ell)$}
\put(2.75,0.3){\ft$\CC$}
\put(-4.3,-1.6){\ft$\Q$}
\put(-0.6,-1.6){\ft$\Q_{n_0}\!=\!\Q_\infty\! \cap\! K^\gamma$}
\put(7.6,-1.6){\ft$K^\gamma$}
\put(-2.95,0.35){\ft$\CC_0 \!=\! k$}
\put(9.75,0.35){\ft$K$}
\put(5.8,-0.23){\ft${\tau_\ell^{}}$}
\bezier{300}(3.2,0.3)(6.3,-0.15)(9.4,0.3)
\put(-4.8,0.0){\ft$p^N$}
\put(2.65,3.8){\ft$M\!=M_0\!$}
\put(9.75,3.8){\ft$L$}
\bezier{300}(-4.0,-1.8)(1.8,-3.0)(7.6,-1.8)
\put(10.0,2.15){\ft$G$}
\put(1.8,-2.8){\ft$\Gamma$}
\put(8.9,-0.7){\ft$\gamma$}
\put(4.6,-1.0){\ft$\CG_K$}
\bezier{300}(-3.6,-1.35)(3.1,-1.2)(9.6,0.3)
\end{picture} }} $$
\vspace{1.1cm}
\unitlength=1.0cm

Under capitulations, we still have the fundamental relation \eqref{fondamental}
for which we have to prove the inequality:
\begin{equation*}
(\CE_K : \CE_{\wh K} \CF_{\!K})_{\varphi_0^{}} 
\geq \order \CH^\ar_{K,\varphi} 
= \frac{\order \wh\BJ \big(\CH_{\wh L,{\varphi_0^{}}}^\ram \big)}
{\order \CH_{\wh L,{\varphi_0^{}}}^\ram} \times
\frac{ (\CE_K : \Norm (\CE_L))_{\varphi_0^{}}}{ (\CE_{\wh K} : \Norm (\CE_{\wh L}))_{\varphi_0^{}}},
\end{equation*}

\noindent
whence  to prove, since the first factor is $\big(\order \Ker(\wh \BJ)_{\varphi_0^{}}\big)^{-1}$, that:
$$\order \Ker(\wh \BJ)_{\varphi_0^{}} \times
\frac{( \CE_{\wh K} : \Norm (\CE_{\wh L}))_{\varphi_0^{}}
(\CE_K : \CE_{\wh K} \CF_{\!K})_{\varphi_0^{}}}
{(\CE_K : \Norm (\CE_L))_{\varphi_0^{}}} \geq 1.$$

We have the exact sequence:
\begin{equation*}
\begin{aligned}
1 \to \CE_{\wh K}\Norm (\CE_L)/\Norm (\CE_L) = 
\CE_{\wh K}/\CE_{\wh K} \cap \Norm (\CE_L)&  \too \\
\CE_K/\Norm (\CE_L) &\too \CE_K/ \CE_{\wh K} \Norm (\CE_L) \to 1
\end{aligned}
\end{equation*}

\noindent
and the inclusions $\Norm (\CE_{\wh L}) \subseteq \CE_{\wh K}\cap \Norm (\CE_L)
\subseteq \CE_{\wh K}$, giving:
\begin{equation*}
\begin{aligned}
\order \Ker(\wh \BJ)_{\varphi_0^{}} & \times
\frac{( \CE_{\wh K} : \Norm (\CE_{\wh L}))_{\varphi_0^{}} 
(\CE_K : \CE_{\wh K} \CF_{\!K})_{\varphi_0^{}}}
{(\CE_K : \Norm (\CE_L))_{\varphi_0^{}}}  \\
= \order \Ker(\wh \BJ)_{\varphi_0^{}} & \times
 \frac{(\CE_{\wh K} : \Norm (\CE_{\wh L}))_{\varphi_0^{}}
(\CE_K : \CE_{\wh K} \CF_{\!K})_{\varphi_0^{}}}
{(\CE_{\wh K} : \CE_{\wh K} \cap \Norm (\CE_L))_{\varphi_0^{}}\, %
(\CE_K : \CE_{\wh K} \Norm(\CE_L) )_{\varphi_0^{}}} \\
 = \order \Ker(\wh \BJ)_{\varphi_0^{}} & \times
 \frac{(\CE_{\wh K}\cap \Norm (\CE_L) : \Norm (\CE_{\wh L}))_{\varphi_0^{}}
(\CE_K : \CE_{\wh K} \CF_{\!K})_{\varphi_0^{}} }
{(\CE_K : \CE_{\wh K} \Norm (\CE_L) )_{\varphi_0^{}}}\, \cdot
\end{aligned}
\end{equation*}

We have $\CF_K^{1-\tau_\ell^{-1}} = \Norm (\CF_L) \subseteq \Norm (\CE_L) $, 
which is not necessarily true for $\CF_K$ and prevents to conclude, contrary 
to the case $\ell$ totally inert. 

\smallskip
If this case of RAMC is true, because of the product formula of Theorem 
\eqref{isotopicphi}, one should obtain, under capitulation of $\CH_K$ in $L$:
\begin{equation}\label{egal1}
\order \Ker(\wh \BJ)_{\varphi_0^{}} \! \times \! 
\frac{(\CE_{\wh K}\cap \Norm (\CE_L) : \Norm (\CE_{\wh L}))_{\varphi_0^{}}
(\CE_K : \CE_{\wh K} \CF_{\!K})_{\varphi_0^{}} }
{(\CE_K : \CE_{\wh K} \Norm (\CE_L) )_{\varphi_0^{}}} = 1.
\end{equation}

Numerical investigations will be necessary to study the behavior of each factor, 
but the relations obtained for the case $\CC_0 \varsubsetneqq k$ (\S\,\ref{known}):

\smallskip
$\ \bullet$ $\Norm (\CE_{\wh L,\varphi_0^{}}) \CF_{\!K,\varphi_0^{}} = \Norm (\CE_{L,\varphi_0^{}})$,

\smallskip
$\ \bullet$ $\Ker(\wh \BJ)_{\varphi_0^{}} = 1$,

\smallskip\noindent
with the supplementary one:

\smallskip
$\ \bullet$ $\CE_{\wh K,{\varphi_0^{}}} \! \cap \Norm (\CE_{L,{\varphi_0^{}}}) =
\Norm (\CE_{\wh L,{\varphi_0^{}}})$,

\smallskip\noindent
satisfies the relation \eqref{egal1}, but this is only speculation.

\end{document}